\documentclass{article}
\usepackage{amsfonts}
\usepackage{amsmath}

\input{tcilatex}

\begin{document}

\title{Extended Weyl Calculus and Application to the Phase-Space Schr\"{o}%
dinger Equation}
\author{Maurice A. de Gosson \\
Universit\"{a}t Potsdam, Inst. f. Mathematik \\
Am Neuen Palais 10, Pf 60 15 33\\
D-14415 Potsdam (Germany)\\
E-mail: maurice.degosson@gmail.com}
\date{}
\maketitle

\begin{abstract}
We show that the Schr\"{o}dinger equation in phase space proposed by
Torres-Vega and Frederick is canonical in the sense that it is a natural
consequence of the extended Weyl calculus obtained by letting the Heisenberg
group act on functions (or half-densities) defined on phase space. This
allows us, in passing, to solve rigorously the TF equation for all quadratic
Hamiltonians.
\end{abstract}

\section{Introduction}

The theory of phase-space Schr\"{o}dinger equation 
\begin{equation}
i\hbar \frac{\partial }{\partial t}\Psi (x,p,t)=H(x+i\hbar \partial
_{p},-i\hbar \partial _{x})\Psi (x,p,t)  \label{TVF}
\end{equation}%
corresponding to the quantization rules%
\begin{equation}
x_{j}\longmapsto x_{j}+i\hbar \frac{\partial }{\partial p_{j}}\text{ },\text{
}p_{j}\longmapsto -i\hbar \frac{\partial }{\partial x_{j}}  \label{QR}
\end{equation}%
proposed by Torres-Vega and Frederick \cite{TV1,TV2,MJ} (TF) has deservedly
received much attention these last years (see for instance \cite%
{CM,JQZ,QJ,QGL,MJ,XQ}). TF arrive at their equation by using a generalized
version of the Husimi transform in which coherent states are used as
pseudo-transition matrix elements via a certain kernel $K_{cs}$. (We take
the opportunity to point out that this kernel has been extensively studied
in the mathematical literature under the names of \textquotedblleft
wavepacket transform\textquotedblright\ or \textquotedblleft FBI\
transform\textquotedblright ; see \textit{e.g. }\cite{Naza}). We intend to
show that TF's choice corresponds to a natural extension of Weyl calculus,
and is therefore in a sense canonical.

\subsection*{Notations}

We will use the collective notations $x=(x_{1},...,x_{n})$, $%
p=(p_{1},...,p_{n})$, and $z=(x,p)$. We denote by $\sigma $ the canonical
symplectic form on the phase space $\mathbb{R}_{z}^{2n}=\mathbb{R}%
_{x}^{n}\times \mathbb{R}_{p}^{n}$:%
\begin{equation*}
\sigma (z,z^{\prime })=(z^{\prime })^{T}Jz\text{ \ \ , \ \ }J=%
\begin{bmatrix}
0 & I \\ 
-I & 0%
\end{bmatrix}%
\text{ }
\end{equation*}%
if $z=(x,p)$, $z^{\prime }=(x^{\prime }p^{\prime })$. We denote by $Sp(n)$
the real symplectic group; it consists of all linear automorphisms $S$ of $%
\mathbb{R}_{z}^{2n}$ such that $S^{T}JS=SJS^{T}=J$. $\mathcal{S}(\mathbb{R}%
^{n})$ is the Schwartz space of rapidly decreasing functions on $\mathbb{R}%
^{n}$ and its dual $\mathcal{S}^{\prime }(\mathbb{R}^{n})$ the space of
tempered distributions. The scalar product of two $n$-vectors $x$ and $p$ is
denoted by $xp$, and we will write $\partial _{z}=(\partial _{x},\partial
_{p})$ where%
\begin{equation*}
\partial _{x}=\left( \tfrac{\partial }{\partial x_{1}},\cdot \cdot \cdot ,%
\tfrac{\partial }{\partial x_{n}}\right) \text{ },\text{ }\partial
_{p}=\left( \tfrac{\partial }{\partial p_{1}},\cdot \cdot \cdot ,\tfrac{%
\partial }{\partial p_{n}}\right) .
\end{equation*}

\section{Heisenberg--Weyl Operators on Phase Space\label{secun}}

Let us recall the definition of the Heisenberg--Weyl operators from the
\textquotedblleft Schr\"{o}dingerian\textquotedblright\ point of view (see 
\textit{e.g.} \cite{Littlejohn}). The flow $(f_{t})$ determined by the
Hamiltonian\ $H_{z_{0}}(z)=\sigma (z,z_{0})$ is given by $f_{t}(z)=z+tz_{0}$
hence the time-one mapping $f_{1}$ is just the phase-space translation $%
T(z_{0}):z\longmapsto z+z_{0}$; the action of that translation on functions $%
\Psi =\Psi (z)$ is defined by \textquotedblleft
push-forward\textquotedblright ;%
\begin{equation*}
T(z_{0})\Psi (z)=\Psi (z-z_{0}).
\end{equation*}%
The quantized version of the Hamiltonian $H_{z_{0}}$ is the operator $%
\widehat{H}_{z_{0}}=\sigma (\widehat{z},z_{0})$ with \ $\widehat{z}%
=(x,-i\hbar \partial _{x})$, and the solution of the corresponding Schr\"{o}%
dinger equation%
\begin{equation*}
i\hbar \frac{\partial \psi }{\partial t}=\sigma (\widehat{z},z_{0})\psi 
\text{ \ , \ }\psi (x,0)=\psi _{0}(x)
\end{equation*}%
is the function 
\begin{equation*}
\psi (x,t)=\exp \left[ \frac{i}{\hbar }(tp_{0}x-\frac{t^{2}}{2}p_{0}x_{0})%
\right] \psi _{0}(x-tx_{0}).
\end{equation*}%
The value of $\psi $ a point $x$ at time $t=1$ is denoted by $\widehat{T}%
(z_{0})\psi (x)$: 
\begin{equation}
\widehat{T}(z_{0})\psi (x)=e^{\frac{i}{\hbar }(p_{0}x-\frac{1}{2}%
p_{0}x_{0)}}\psi (x-x_{0.});  \label{defusual}
\end{equation}%
$\widehat{T}(z_{0})$ is by definition the \textit{Heisenberg--Weyl} operator
associated to the translation $T(z_{0})$ (see \cite{Littlejohn}). So much
for the traditional point of view. Let us shortly review the notion of phase
as studied in our paper \cite{phases}. Consider a Lagrangian manifold $%
\mathbb{V}_{0}^{n}$ in $\mathbb{R}_{z}^{2n};$ we assume for simplicity that $%
\mathbb{V}_{0}^{n}$ is simply connected (this restriction, which can be
alleviated by passing to the universal covering, is of no consequence for
the argument that follows). Choosing, once for all, a base point $\bar{z}$
in $\mathbb{V}_{0}^{n}$ we define the \textit{phase} of $\mathbb{V}_{0}^{n}$
by%
\begin{equation*}
\varphi _{0}(z)=\int_{\gamma (\bar{z},z)}pdx\text{ }
\end{equation*}%
where the integral is calculated along an arbitrary path leading from $\bar{z%
}$ to $z.$ This defines $\varphi _{0}$ as a function $\mathbb{V}%
_{0}^{n}\longrightarrow \mathbb{R}$ such that $d\varphi _{0}(z)=pdx$\ if $%
z=(x,p)$ ($\mathbb{V}_{0}^{n}$ being Lagrangian, the integral only depends
on the homotopy class with fixed endpoints of $\gamma (\bar{z},z)$ and since 
$\mathbb{V}_{0}^{n}$ is in addition simply connected there is just one
homotopy class). Let now $H$ be some arbitrary Hamiltonian function; its
flow $(f_{t})$ takes $\mathbb{V}_{0}^{n}$ to a new Lagrangian manifold $%
\mathbb{V}_{t}^{n}$ after time $t$. The phase of that manifold is (with
obvious notations)%
\begin{equation*}
\varphi (z,t)=\varphi _{0}(f_{-t}(z))+\int_{\Gamma }pdx-Hdt
\end{equation*}%
where the integral is calculated along the arc of phase space trajectory
leading from the point $f_{-t}(z)\in \mathbb{V}_{0}^{n}$ at time $0$ to the
point $z\in \mathbb{V}_{t}^{n}$ at time $t$. Choose now for $H=H_{z_{0}}$; a
straightforward calculation (see \cite{phases}, Prop. 10) yields%
\begin{equation*}
\varphi (z,t)=\varphi _{0}(z-z_{0})+tp_{0}x-\frac{t^{2}}{2}p_{0}x_{0}.
\end{equation*}%
We next consider, as is customary in geometric quantization (see \cite%
{Bullsci,IHP,ICP}), a \textquotedblleft waveform\textquotedblright\ on $%
\mathbb{V}_{0}^{n}$, \textit{i.e.} an expression of the type%
\begin{equation*}
\Psi _{0}(z)=e^{\frac{i}{\hbar }\varphi _{0}(z)}i^{m_{0}(z)}\sqrt{\rho }(z)
\end{equation*}%
where $m_{0}(z)$ is an integer related to the Maslov index and $\rho $ a
function $\geq 0$ on$\mathbb{V}_{0}^{n}$ (more accurately, one should view $%
\sqrt{\rho }$ as a \textquotedblleft half-density\textquotedblright\ on $%
\mathbb{V}_{0}^{n}$ but we can ignore this point here because half-densities
and functions transform the same way under translations). The action of a
Hamiltonian flow $(f_{t})$ on $\Psi _{0}$ is given by $f_{t}\Psi
_{0}(z)=\Psi (z,t)$ where%
\begin{equation*}
\Psi (z,t)=e^{\frac{i}{\hbar }\varphi (z,t)}i^{m(z,t)}(f_{t})_{\ast }\sqrt{%
\rho }(z)
\end{equation*}%
the difference $m(z,t)-m_{0}(z)$ is a count of the number of caustic points
appearing when the flow deforms $\mathbb{V}_{0}^{n}$ into $\mathbb{V}%
_{t}^{n}=f_{t}(\mathbb{V}_{0}^{n})$. If $H=H_{z_{0}}$ the manifold $\mathbb{V%
}_{0}^{n}$ is translated in phase space without any \textquotedblleft
bending\textquotedblright\ so that there will appear no new caustics hence $%
m(z,t)=m_{0}(z)$ and we thus have%
\begin{equation}
\Psi (z,t)=\hat{T}(tz_{0})\Psi _{0}(z)  \label{tete}
\end{equation}%
where the operator $\hat{T}(z_{0})$ is now defined by%
\begin{equation}
\hat{T}(z_{0})\Psi _{0}(z)=e^{\frac{i}{\hbar }(p_{0}x-\frac{1}{2}%
p_{0}x_{0})}T(z_{0})\Psi _{0}(z).  \label{WHP}
\end{equation}

A straightforward calculation shows that the partial differential equation
satisfied by $\Psi (z,t)$ is 
\begin{equation}
i\hbar \frac{\partial \Psi }{\partial t}(z,t)=\widehat{\Sigma }\Psi (z,t)%
\text{ },\text{ }\Psi (t=0)=\Psi _{0}  \label{popo}
\end{equation}%
where $\widehat{\Sigma }$ is the operator $\widehat{\Sigma }=-p_{0}x+i\hbar
z_{0}\partial _{z};$the latter can be rewritten in terms of the symplectic
form $\sigma $ as%
\begin{equation*}
\widehat{\Sigma }=\sigma (x_{0},p_{0};x+i\hbar \partial _{p},-i\hbar
\partial _{x})
\end{equation*}%
so that (\ref{popo}) is just TF's phase-space Schr\"{o}dinger equation (\ref%
{TVF}) for the Hamiltonian $H_{z_{0}}$.

The whole argument can actually be reversed: using for instance the methods
of characteristics one checks that the solution of equation (\ref{popo}) is
precisely (\ref{tete}) so that we conclude that the imposition of the TF
quantization rules (\ref{QR}) is \emph{equivalent} to the extension (\ref%
{WHP}) of Heisenberg--Weyl operators to phase-space.

\section{Weyl Calculus on Phase Space\label{secdeux}}

The process called \textquotedblleft Weyl quantization\textquotedblright\
associates to a suitable function (or \textquotedblleft
symbol\textquotedblright ) $a=a(z)$ the operator $a^{w}=\hat{A}$ defined by%
\begin{equation}
\hat{A}\psi (x)=\left( \tfrac{1}{2\pi \hbar }\right) ^{n}\dint \tilde{a}%
(z_{0})\hat{T}(z_{0})\psi (x)d^{2n}z_{0}  \label{weyl1}
\end{equation}%
where $\tilde{a}$ is the \textquotedblleft alternative Weyl
symbol\textquotedblright , that is the symplectic Fourier transform $%
\mathcal{F}_{\sigma }a$ of $a$:%
\begin{equation*}
\tilde{a}(z)=\mathcal{F}_{\sigma }a(z)=\left( \tfrac{1}{2\pi \hbar }\right)
^{n}\dint e^{-\frac{i}{\hbar }\sigma (z,z^{\prime })}a(z^{\prime
})d^{2n}z^{\prime }.
\end{equation*}

The discussion of Section \ref{secun} suggests that we might now be able to
make $\hat{A}$ to act, not only on function of $x,$ but also on functions of 
$z$ by defining 
\begin{equation}
\hat{A}\Psi (z)=\left( \tfrac{1}{2\pi \hbar }\right) ^{n}\dint \tilde{a}%
(z_{0})\widehat{T}(z_{0})\Psi (z)d^{2n}z_{0}  \label{weyl2}
\end{equation}%
where $\Psi $ now is a function of $z=(x,p)$ and $\widehat{T}(z_{0})\Psi (z)$
is given by (\ref{WHP}); this expression makes perfectly sense provided of
course, that $a$ and $\Psi $ belong to some adequate spaces; we will assume
this is the case. It turns out that this reinterpretation of Weyl operators
again leads to the TF quantization rules (\ref{QR}). It is sufficient to
prove this in the case of one degree of freedom; dropping indices let $%
\widehat{X}$ and $\widehat{P}$ be the Weyl operators with symbols $x$ and $p$%
; a simple calculation of (symplectic) Fourier transforms yields%
\begin{eqnarray*}
\widehat{X}\Psi (z) &=&i\hbar \int \delta (x_{0})\otimes \delta ^{\prime
}(p_{0})\widehat{T}(z_{0})\Psi (z)d^{2n}z_{0} \\
&=&-i\hbar \int \delta (z_{0})(\tfrac{ix}{\hbar }\Psi (z)-\tfrac{\partial
\Psi }{\partial p}(z))d^{2n}z_{0} \\
&=&(x+i\hbar \tfrac{\partial }{\partial p})\Psi (z)
\end{eqnarray*}%
and we thus have $\widehat{X}=x+i\hbar \tfrac{\partial }{\partial p};$ by a
similar calculation we get $\widehat{P}=-i\hbar \frac{\partial }{\partial x}$%
, and we have thus recovered the TF rules (\ref{QR}).

\section{Exact Solutions for Quadratic Hamiltonians}

The extended Weyl calculus constructed in previous Section allows us to
solve exactly the TF equation when the Hamiltonian $H$ is a homogeneous
quadratic polynomial in the variables $x_{j},p_{k}$, in which case the
Hamiltonian flow consists of symplectic matrices $S_{t}$.

One of the main features (and probably one of the most attractive!) of the
usual Weyl correspondence $a\longleftrightarrow a^{w}=\hat{A}$ is the
property of \emph{metaplectic covariance}: if we replace $a$ by $a\circ
S^{-1}$ where $S$ is any element of $Sp(n)$ then :%
\begin{equation}
(a\circ S^{-1})^{w}=\hat{S}a^{w}\hat{S}^{-1}  \label{momo}
\end{equation}%
where $\hat{S}$ is any of the two metaplectic operators $\pm \hat{S}$
associated to $S$ (see \cite{Littlejohn}, Section 6.3; for the properties of 
$Mp(n)$ see \cite{ICP}, Chapter 6). The proof of property (\ref{momo}) is
based on the formula%
\begin{equation}
\hat{S}\hat{T}(z_{0})\hat{S}^{-1}=\hat{T}(Sz_{0})  \label{chacha}
\end{equation}%
where one uses the traditional definition (\ref{defusual}) definition of $%
\hat{T}(z_{0})$, hence if (\ref{chacha}) still holds when we redefine $\hat{T%
}(z_{0})$ by (\ref{WHP}), then metaplectic covariance (\ref{momo}) will also
hold for our enlarged Weyl calculus. Now, Mehlig and Wilkinson have shown 
\cite{MW} that if $\det (S-I)\neq 0$ then%
\begin{equation}
\hat{S}=\pm k\int \hat{T}(Sz_{0})\hat{T}(-z_{0})d^{2n}z_{0}  \label{toto}
\end{equation}%
where $k$ is a constant given by%
\begin{equation}
k=\left( \tfrac{1}{2\pi \hbar }\right) ^{n}i^{\nu }\sqrt{\left\vert \det
(S-I)\right\vert }  \label{totobis}
\end{equation}%
and $\nu $ is a kind of \textquotedblleft Maslov index\textquotedblright .
We have shown \cite{MdGMW} that every $\hat{S}\in Mp(n)$ can be written as
the product of two operators of the type (\ref{toto}) above (which are only
defined for $\det (S-I)\neq 0);$ we may thus extend every $\hat{S}\in Mp(n)$
to an operator defined on functions of $z=(x,p);$ the natural domain of $%
\hat{S}$ is then $\mathcal{S}(\mathbb{R}_{z}^{2n}),$ and can be extended to $%
L^{2}(\mathbb{R}^{2n})$ and to $\mathcal{S}^{\prime }(\mathbb{R}_{z}^{2n})$
by continuity. Having done this it is immediate to check, using formulae (%
\ref{toto}) that (\ref{chacha}), and hence also (\ref{momo}), remain true in
our extended Weyl calculus.

Let us now return to the flow $(S_{t})$; denote by $(\hat{S}_{t})$ the
unique one-parameter group in $Mp(n)$ determined by $(S_{t})$; viewing the $%
\hat{S}_{t}$ as extended Weyl operators we automatically have 
\begin{equation*}
i\hbar \frac{d}{dt}\hat{S}_{t}=H(x+i\hbar \partial _{p},-i\hbar \partial
_{x})\hat{S}_{t}
\end{equation*}%
hence $\Psi (z,t)=\hat{S}_{t}\Psi _{0}(z)$ is the solution of the TF
equation (\ref{TVF}) for the considered Hamiltonian.

\section{Discussion and Concluding Remarks}

We have established that the Torres-Vega and Frederick's equation fits as
part of an extended Weyl calculus, which was itself motivated by
considerations from geometric quantization. There remains to find the \emph{%
physical interpretation} of that equation. A clue might be the fact that it
is possible to rewrite (a variant of) the TF\ equation in terms of the
star-product, familiar from deformation quantization (this fact seems to
have been noticed in \cite{CM}; this and the fact that the star-product is
used to express the time-evolution of the Wigner function via a
quantum-Liouville equation might be an indication that the TF equation plays
in quantum mechanics a role similar to that of Hamilton's equation in
classical mechanics.

We will come back to this important question in a near future.


\begin{thebibliography}{99}
\bibitem{CM} D. Chruscinski and K. Mlodawski. Preprint, arXiv: quant-ph/
0501163 v1, 2005.

\bibitem{Bullsci} M. de Gosson. \textit{Bull. Sci. Math}. \textbf{121}%
:301--322, 1997.

\bibitem{IHP} M. de Gosson. \textit{Ann. Inst. H. Poincar\'{e}}, \textbf{70}%
(6):547--73, 1999.

\bibitem{ICP} M. de Gosson. \textit{The Principles of Newtonian and Quantum
Mechanics. }Imperial College Press, London, 2001.

\bibitem{phases} M. de Gosson. \textit{J. Phys. A: Math. Gen}. \textbf{37}%
(29), 7297--7314, 2004.

\bibitem{MdGMW} M. de Gosson. On the Weyl Representation of Metaplectic
Operators. To appear in \textit{Lett. Math. Phys.,} 2005

\bibitem{JQZ} Jun Lu, Qiam-Shu Li and Zheng-Sun. \textit{J. Chem. Phys.} 
\textbf{3}, 1022--1036, 2001.

\bibitem{Littlejohn} R. G. Littlejohn. \textit{Physics Reports} \textbf{138(}%
4--5):193--291, 1986.

\bibitem{QJ} Qiam-Shu Li, Jun Lu. \textit{Chem. Phys. Letter}. \textbf{336},
118--122, 2001.

\bibitem{QGL} Qiam-Shu Li, Gong Min Wei, and Li Qiang Lu. \textit{Phys. Rev.}%
\textbf{\ }\textit{A}. \textit{70}, 022105 (1--5), 2004.

\bibitem{MW} B. Mehlig and M. Wilkinson\textsc{.} \textit{Ann. Phys.} 
\textbf{18}(10), 6--7, 541--555, 2001.

\bibitem{MJ} K. B. Moeller, T. G. Joergensen and G. Torres-Vega. \textit{J.
Chem. Phys.} \textbf{106}(17), 7228--7240, 1997.

\bibitem{Naza} V. Nazaikiinskii, B.-W. Schulze, and B. Sternin. \textit{%
Quantization Methods in Differential Equations}. Taylor \& Francis, 2002.

\bibitem{TV1} G. Torres-Vega and J. H. Frederick. \textit{J. Chem. Phys}. 
\textbf{93}(12), 8862--8874, 1990.

\bibitem{TV2} G. Torres-Vega and J. H. Frederick. \textit{J.} \textit{Chem.
Phys.} \textbf{98}(4), 3103--3120, 1993.

\bibitem{XQ} Xu-Guang Hu and Qiam-Shu Li. \textit{J. Phys. A: Math. Gen. }%
\textbf{32}, 139--146, 1999.
\end{thebibliography}
\end{document}